\documentclass[12pt,a4paper]{article}

\usepackage{amsmath,amssymb}

\def\cB{{\cal B}}

\def\cF{{\cal F}}
\def\cK{{\cal K}}
\def\cS{{\cal S}}

\def\Z{{\mathbb Z}}

\def\C{{\mathbb C}}

\def\F{{\mathbb F}}
\def\sl{{\mathfrak sl}}

\def\corank{{\rm corank}}

\usepackage{theorem}

\newtheorem{theorem}{Theorem}[section]
\newtheorem{proposition}{Proposition}[section]
\newtheorem{corollary}[proposition]{Corollary}

{\theorembodyfont{\rmfamily}
\newtheorem{definition}[proposition]{Definition}
\newtheorem{example}[proposition]{Example}
\newtheorem{remark}[proposition]{Remark}

}

\usepackage{graphicx}

\title{Delta-matroids and Vassiliev invariants}

\author{
Sergey Lando\thanks{National Research University Higher School of Economics,
Independent University of Moscow, lando@hse.ru},
Vyacheslav Zhukov\thanks{National Research University Higher School of Economics, slava.zhukov@list.ru}}

\date{}

\begin{document}
\maketitle

Vassiliev (finite type) invariants of knots can be described in terms of weight systems.
These are functions on chord diagrams satisfying so-called $4$-term relations.
In the study of the $\sl_2$ weight system in~\cite{CL07}, it was shown that
its value on a chord diagram depends on the intersection graph of the diagram rather
than on the diagram itself. Moreover, it was shown that the value of this
weight system on an intersection graph depends on the cycle matroid of the graph
rather than on the graph itself. This result arose the question
whether there is a natural way to introduce a $4$-term relation on the space spanned
by matroids, similar to the one for graphs~\cite{L00}. It happened however that the
answer is negative: there are graphs having isomorphic cycle matroids such that
applying the ``second Vassiliev move'' to a pair of corresponding vertices
$a,b$ of the graphs we obtain two graphs with nonisomorphic matroids.

The goal of the present paper is to show that the situation is different for
binary delta-matroids: one can define both the first and the second Vassiliev
moves for binary delta-matroids and introduce a $4$-term relation for them in such
a way that the mapping taking a chord diagram to its delta-matroid
respects the corresponding $4$-term relations. Moreover, this mapping admits
a natural extension to chord diagrams on several circles, which
correspond to singular links. Delta-matroids were introduced by A.~Bouch\'et~\cite{Bo871}
for the purpose of studying embedded graphs, whence their relationship
with (multiloop) chord diagrams is by no means unexpected. Some evidence
for the existence of such a relationship can be found, for example, in~\cite{BR01},
where the Tutte polynomial for embedded
graphs has been introduced. The authors show that this polynomial depends on the
delta-matroid of the embedded graph rather than the graph itself and satisfies the Vassilev $4$-term relation.

Understanding how the $4$-term relation can be written out for
arbitrary binary delta-matroids motivates introduction of the graded Hopf algebra of binary delta-matroids
modulo the $4$-term relations so that the mapping taking a chord diagram
to its delta-matroid extends to a morphism of Hopf algebras.
One can hope that studying this Hopf algebra will allow one to clarify
the structure of the Hopf algebra of weight systems, in particular, to find
reasonable new estimates for the dimensions of the spaces of weight systems
of given degree. Also it would be interesting to find a relationship between the Hopf algebras
arising in this paper with a very close to them in spirit bialgebra of
Lagrangian subspaces in~\cite{KS}.

The authors are grateful to participants of the seminar ``Combinatorics of Vassiliev invariants''
at the Department of mathematics, Higher School of Economics and Sergei Chmutov for
useful discussions. The article was prepared within the framework of the Academic Fund Program
at the National Research University Higher School of Economics (HSE) in 2016---2016 (grant ¹ 16-05-0007)
and supported within the framework of a subsidy granted to the HSE by the Government
of the Russian Federation for the implementation of the Global Competitiveness Program.

\section{Algebra of set systems}\label{sHa}

A {\it set system\/} $(E;\Phi)$ is a finite set~$E$ together with a subset~$\Phi$
of the set $2^E$ of subsets in~$E$. The set~$E$ is called the {\it ground
set\/} of the set system, and elements of~$\Phi$ are its {\it feasible sets}.
Two set systems $(E_1;\Phi_1)$, $(E_2;\Phi_2)$ are said to
be {\it isomorphic\/} if there is a one-to-one map $E_1\to E_2$ identifying
the subset~$\Phi_1\subset 2^{E_1}$ with the subset~$\Phi_2\subset 2^{E_2}$. Below, we make no difference between isomorphic set systems.

A set system $(E;\Phi)$ is {\it proper\/} if~$\Phi$ is nonempty. Below, we consider only
proper set systems, without indicating this explicitly.

\subsection{The graded vector space of set systems}

Let $\cS_n$ denote the vector space (over the field of complex numbers~$\C$, for definiteness)
freely spanned by set systems whose ground set consists of~$n$ elements,
$\cS_0$ being the field~$\C$ itself.
The direct sum
$$
\cS=\cS_0\oplus \cS_1\oplus \cS_2\oplus\dots
$$
is an infinite dimensional graded vector space.

\begin{example}
The vector space~$\cS_0$ is $1$-dimensional. It is spanned by the only set system
on zero elements, namely, the set system $\{\emptyset;\{\emptyset\}\}$.

The vector space~$\cS_1$ is $3$-dimensional. It is spanned by the three set systems
$$
s_{11}=\{\{1\};\{\emptyset\}\},\quad
s_{12}=\{\{1\};\{\emptyset,\{1\}\}\},\quad
s_{13}=\{\{1\};\{\{1\}\}\}.
$$
\end{example}

In our notation~$s_{ij}$ for set systems, the first index~$i$ denotes the number of
elements in the ground set, while the second one is chosen ambiguously.

\begin{remark}
Note that the set systems $\{\emptyset;\{\emptyset\}\}$
and $s_{11}$ {\it are\/} proper. Indeed, in both cases the corresponding set of
subsets is not empty: it contains one element, namely, the empty set.

\end{remark}

\subsection{Multiplication of set systems}

The {\it direct sum\/} of two set systems $D_1=(E_1;\Phi_1)$, $D_2=(E_2;\Phi_2)$
with disjoint ground sets $E_1,E_2$ is defined to be
\begin{equation}\label{emult}
D_1D_2=(E_1\sqcup E_2;\{\phi_1\sqcup \phi_2|\phi_1\in \Phi_1,\phi_2\in\Phi_2\}).
\end{equation}
Since we consider set systems up to isomorphism, we will always assume that,
when considering direct sums, the ground sets $E_1$ and $E_2$ of the summands are disjoint.
Below, we will also refer to the direct sum as to the {\it product\/}
of set systems. This operation extends by linearity to a bilinear multiplication
$$
m:\cS\otimes\cS\to\cS,\qquad m(D_1\otimes D_2)=D_1D_2,
$$
which is graded (meaning that $m:\cS_k\otimes \cS_\ell\to \cS_{k+\ell}$ for all $k,\ell\ge0$),
and commutative. The unit of this multiplication is the set
system $(\emptyset;\{\emptyset\})$, which is the generator of~$S_0$.

\begin{example}
The vector space~$\cS_2$ is $11$-dimensional. It is spanned by the six set systems
that are products of set systems on one element sets, namely,
\begin{eqnarray*}
s_{11}^2&=&\{\{1,2\};\{\{\emptyset\}\}\},\\
s_{12}^2&=&\{\{1,2\};\{\emptyset,\{1\},\{2\},\{1,2\}\}\},\\
s_{13}^2&=&\{\{1,2\};\{\{1,2\}\}\},\\
s_{11}s_{12}&=&\{\{1,2\};\{\emptyset,\{1\}\}\}=\{\{1,2\};\{\emptyset,\{2\}\}\},\\
s_{11}s_{13}&=&\{\{1,2\};\{\{1\}\}\}=\{\{1,2\};\{\{2\}\}\},\\
s_{12}s_{13}&=&\{\{1,2\};\{\{1\},\{1,2\}\}\}=\{\{1,2\};\{\{2\},\{1,2\}\}\},
\end{eqnarray*}
and the five other set systems
\begin{eqnarray*}
s_{21}&=&\{\{1,2\};\{\emptyset,\{1,2\}\}\},\\
s_{22}&=&\{\{1,2\};\{\emptyset,\{1\},\{1,2\}\}\}=\{\{1,2\},\{\emptyset,\{2\},\{1,2\}\}\},\\
s_{23}&=&\{\{1,2\};\{\emptyset,\{1\},\{2\}\}\},\\
s_{24}&=&\{\{1,2\};\{\{1\},\{2\}\}\},\\
s_{25}&=&\{\{1,2\};\{\{1\},\{2\},\{1,2\}\}\}.
\end{eqnarray*}

\end{example}

\section{Generalities on delta-matroids}

In this section we briefly reproduce the general facts about delta-matroids that we will
require further. We follow the approach and terminology of~\cite{CMNR}, but use
slightly different notation.

\subsection{Delta-matroids}

Let~$\Delta$ denote the symmetric difference of sets, $A\Delta B=(A\setminus B)\sqcup (B\setminus A)$.
A {\it delta-matroid\/} is a set system $D=(E;\Phi)$ satisfying the following
Symmetric Exchange Axiom (SEA):

{\it For any $\phi_1,\phi_2\in \Phi$ and for any $e\in(\phi_1\Delta\phi_2)$ there is
an element $e'\in(\phi_2\Delta\phi_1)$ such that $\phi_1\Delta\{e,e'\}\in\Phi$.
}


It is easy to check that all the set systems on~$1$ or~$2$ elements, which are enumerated
in Sec.~\ref{sHa}, are delta-matroids. However, there are set systems
that are not delta-matroids already among set systems on three elements.
For example, if, for the set system
$(\{1,2,3\};\{\emptyset,\{1,2,3\}\})$ we take $\phi_1=\emptyset$,
$\phi_2=\{1,2,3\}$, then the SEA will not be satisfied.

\subsection{Delta-matroids of embedded graphs}
An {\it embedded graph\/} is, essentially, a graph drawn on a
compact surface in such a way that its complement is a disjoint
union of disks. We will always assume that the graph is connected.
Edges in an embedded graph are also called {\it ribbons}, or {\it handles},
and we make no distinction between embedded and ribbon graphs.
Generalities on embedded graphs can be found, for example, in~\cite{LZ04}.

If otherwise is not stated explicitly, then we allow both orientable and
nonorientable surfaces. A loop in an embedded graph, that is, an edge connecting
a vertex with itself, can be orientable or disorienting (half-twisted).
If there is a disorienting loop in an embedded graph, then the graph itself
is nonorientable. However, a nonorientable ribbon graph does not necessarily
contain a disorienting loop: it suffices that there exists a disorienting cycle,
not necessarily of length~$1$, in it.

To each embedded graph~$\Gamma$, its delta-matroid $D(\Gamma)=(E(\Gamma);\Phi(\Gamma))$ is associated.
The ground set of the delta-matroid is the set $E(\Gamma)$ of the edges of~$\Gamma$.
A subset $\phi\subset E(\Gamma)$ is feasible, $\phi\in\Phi(\Gamma)$,
 if the boundary of the embedded spanning subgraph
of~$\Gamma$ formed by the set~$\phi$ is connected, that is, consists of a single connected component.
This means, in particular, that the spanning subgraph of~$\Gamma$
formed by the set~$\phi$ is connected (otherwise, each connected component
would add at least one connected component to the boundary). Since, for a plane graph, this
requirement coincides with the requirement that~$\phi$ is a spanning tree,
feasible sets for graphs embedded into a surface of arbitrary genus are called
{\it quasi-trees}. For graphs embedded in surfaces of positive genus, not all
of quasi-trees necessarily are trees, although each subset of edges forming
a spanning tree is feasible.

Delta-matroids of orientable embedded graphs are {\it even}, meaning that all
the feasible sets in them have cardinality of the same parity.

All the set systems in Sec.~\ref{sHa} are delta-matroids of embedded graphs.
Thus, $s_{11}$ is the delta-matroid of the embedded graph with one vertex
and an orientable loop,  $s_{12}$ is the delta-matroid of the embedded graph with one vertex
and a half-twisted loop, while $s_{13}$ is the delta-matroid of the embedded graph with two vertices
and an edge connecting them. The delta-matroids $s_{11}, s_{13}$ correspond to
orientable embedded graphs, and are even, while~$s_{12}$ is not even.

The following statement is straightforward.

\begin{proposition}[\cite{CMNR}]
If $\Gamma_1$, $\Gamma_2$ are two embedded graphs with the delta-matroids $D(\Gamma_1),D(\Gamma_2)$,
respectively, then the delta-matroid of the embedded graph $\Gamma_1\#\Gamma_2$
obtained by gluing $\Gamma_1,\Gamma_2$ along a vertex is the product of the delta-matroids of
the summands, $D(\Gamma_1\#\Gamma_2)=D(\Gamma_1)D(\Gamma_2)$.
\end{proposition}

Here the gluing $\Gamma_1\#\Gamma_2$ of embedded graphs~$\Gamma_1, \Gamma_2$
along a vertex is defined in the following way:
we choose an arbitrary vertex in~$\Gamma_1$ and an arbitrary vertex in~$\Gamma_2$,
and glue the two vertices together so that the half-edges of~$\Gamma_1$
leave the joint vertex in the same cyclic order, followed by the those of~$\Gamma_2$.
The above proposition means, in particular, that the delta-matroid of the resulting graph depends neither
on the choice of the two vertices to be glued, nor on the choice of the
breaking point inside each vertex. Note that the number of vertices in the result
of gluing of two graphs is one less than the total number of vertices in the graphs.

\begin{example}
The delta-matroid $s_{13}^2$ is represented by the only embedded graph
with three vertices and two edges.
\end{example}

\section{$\Delta$-matroids of abstract graphs and binary delta-matroids}

Certain abstract graphs can be represented as intersection graphs of {\it chord diagrams},
which are embedded graphs with a single vertex. In spite of the fact that
one graph can be the intersection graph of different chord diagrams,
all these diagrams have have one and the same delta-matroid, which is, therefore,
associated to the graph itself.
Bouch\'et extended this construction to arbitrary abstract
graphs.

\subsection{Binary delta-matroids}

Let~$G$ be an (abstract) undirected graph. We say that~$G$ is nondegenerate if
its adjacency matrix $A(G)$, considered as a matrix over the field of two elements,
is nondegenerate. Define the set system $(V(G);\Phi(G))$, $\Phi(G)\subset 2^{V(G)}$,
by
\begin{eqnarray*}
V(G) &\text{ is }&\text{ the set of vertices of } G, \\
\Phi(G)&=&\{U\subset V(G)|G_U \text{ is nondegenerate} \},
\end{eqnarray*}
where $G_U$ is the subgraph in~$G$ induced by the subset~$U$ of vertices.

\begin{theorem}[\cite{Bo871}]\label{tDG}
The set system $(V(G);\Phi(G))$ is a delta-matroid.
\end{theorem}

We call this delta-matroid the {\it nondegeneracy delta-matroid \/} of the graph~$G$.

For an orientable embedded graph~$\Gamma$ with a single vertex, denote by $\gamma(\Gamma)$
its intersection graph, that is, the graph whose vertices correspond one-to one to the
ribbons of~$\Gamma$, and two vertices are connected by an edge iff the ends of the
corresponding ribbons alternate along the vertex.

\begin{theorem}[\cite{Bo871}]\label{tDMIG}
Let~$\Gamma$ be an orientable ribbon graph with a single vertex. Then its $\Delta$-matroid
$(E(\Gamma);\Phi(\Gamma))$ coincides with the nondegeneracy delta-matroid of the intersection
graph $\gamma(\Gamma)$ of~$G$.
\end{theorem}

According to the theorem from~\cite{S01}, the number of connected components of the
boundary of a ribbon graph~$\Gamma$ with a single vertex is equal to $\corank(A(\gamma(\Gamma)))+1$,
where the adjacency matrix is considered over the field with two elements.
In particular, the boundary has a single component iff the matrix $A(\gamma(\Gamma))$ is
nondegenerate.

Theorem~\ref{tDG} is naturally generalized to framed graphs and nonorientable embedded graphs.
Recall the definition of a framed graph from~\cite{L}.

\begin{definition}
A {\it framed graph\/} is an (abstract) graph~$G$ together with a {\it framing}, that is, a mapping
$f:V(G)\to\{0,1\}$. In the adjacency matrix $A(G)$ of a framed graph, the
diagonal element corresponding to a vertex $v\in V(G)$ is $f(v)$, while
nondiagonal elements are defined as usual.
\end{definition}

For a framed graph~$G$, the set system $(V(G);\Phi(G))$, is defined in
the same way as for an unframed one.

Now let~$\Gamma$ be a ribbon graph with a single vertex, not necessarily orientable.
The intersection graph~$\gamma(\Gamma)$ of the ribbon graph~$\Gamma$ is the framed graph
such that each nonoriented loop is taken to a vertex with framing~$1$.
The theorem from~\cite{S01} has the following framed analogue.

\begin{theorem}
For an embedded graph~$\Gamma$ with a single vertex, not necessarily orientable, let
$A(\gamma(\Gamma))$ be the adjacency matrix of its framed intersection graph. Then the
number of connected components of the boundary of~$\Gamma$ is equal to $\corank(A(\gamma(\Gamma)))+1$.
\end{theorem}

As a consequence, we obtain a generalization of Theorem~\ref{tDMIG} for not necessarily
orientable ribbon graph with a single vertex.

\begin{corollary}
Let~$\Gamma$ be a ribbon graph with a single vertex. Then its delta-matroid
$(E(\Gamma);\Phi(\Gamma))$ coincides with the nondegeneracy $\Delta$-matroid of the intersection
graph $\gamma(\Gamma)$ of~$\Gamma$.
\end{corollary}

Nondegeneracy delta-matroids of abstract framed graphs are examples of binary delta-matroids.
In order to define the notion of binary delta-matroid, we will need the twist operation.
For a set system $D=(E;\Phi)$ and a subset $E'\subset E$ define the {\it twist\/} $D*E'$ of~$D$
with respect to~$E'$ by
$$
D*E'=(E;\Phi\Delta E')=(E;\{\phi\Delta E'|\phi\in\Phi\}).
$$

\begin{theorem}[\cite{Bo89}]
Any twist of a nondegeneracy delta-matroid of a framed graph is a delta-matroid.
\end{theorem}

Bouch\'et calls the delta-matroids obtained as twists of nondegenracy delta-matroids
of framed graphs {\it binary delta-matroids\/}. In particular, he shows that

\begin{theorem}[\cite{Bo89}]
Delta-matroids of embedded graphs are binary.
\end{theorem}

Below, we will consider the algebra of binary delta-matroids. It is well-defined due to the
following statement.

\begin{theorem}[\cite{CMNR}]
The product of two binary delta-matroids is a binary delta-matroid.
\end{theorem}

This theorem means that we can consider the graded commutative algebra of binary delta-matroids,
which is a graded subalgebra in the algebra~$\cS$ of set systems. We will denote this algebra by~$\cB$:
$$
\cB=\cB_0\oplus\cB_1\oplus\cB_2\oplus\dots.
$$
The graded subalgebra $\cB^e$ in~$\cB$ is spanned by even binary delta-matroids. Recall
that a delta-matroid $(E;\Phi)$ is {\it even\/} if the parity of the cardinality
is the same for all sets in~$\Phi$.

\subsection{Comultiplication of binary delta-matroids}

In addition to multiplication, we are going to introduce a comultiplication~$\mu$ on the
space~$\cB$ of binary delta-matroids, $\mu:\cB\to\cB\otimes\cB$.
By definition, the {\it coproduct\/} $\mu(D)$ of a delta-matroid $D=(E;\Phi)$ is
\begin{equation}\label{ecomult}
\mu(D)=\sum_{E'\subset E}D_{E'}\otimes D_{E\setminus E'}.
\end{equation}
Here, for a subset $E'\subset E$ of the ground set~$E$ of a delta-matroid~$D$, we denote
by $D_{E'}$ the restriction of~$D$ to~$E'$.

Let us recall the definition of restriction from~\cite{CMNR}. It requires some other notions,
which we collect together in a single paragraph.

\begin{definition}\rm
Let $D=(E;\Phi)$ be a delta-matroid.
An element  $e\in E$ is a {\it coloop\/} if it enters all feasible sets in~$D$.
If~$e$ is not a coloop, then the delta-matroid~$D$ {\it delete\/}~$e$, $D\setminus\{e\}$
is
$$
D\setminus\{e\}=(E\setminus\{e\};\{\phi\in\Phi|\phi\subset E\setminus\{e\}\}).
$$
An element  $e\in E$ is a {\it loop\/} if it does not enter any feasible set in~$D$.
If~$e$ is not a loop, then the delta-matroid~$D$ {\it contract\/}~$e$, $D/\{e\}$
is
$$
D/\{e\}=(E\setminus\{e\};\{\phi\setminus\{e\}|\phi\in\Phi\text{ and }\phi\ni e\}).
$$
If~$e$ is a coloop, then, by definition,  $D\setminus\{e\}=D/\{e\}$.
If~$e$ is a loop, then, by definition,  $D/\{e\}=D\setminus\{e\}$.
A {\it minor\/} of~$D$ is a delta-matroid obtained from~$D$ by a sequence
of deletions and contractions. The {\it restriction\/} $D_{E'}$ of~$D$ to a subset $E'\subset E$
is the result of deleting all elements in $(E\setminus E')\subset E$ in~$D$.
\end{definition}

All these notions are well-defined. This means, in particular, that the deletion and contraction of a
delta-matroid are delta-matroids as well, and any sequence of deletions and contractions leads
to the same delta-matroid independently of the order of the elements in the sequence
(which are assumed to be pairwise distinct). In the notation below, we will often omit braces
around one-element sets, writing $E\setminus e$ instead of $E\setminus\{e\}$, and so on.

\begin{proposition}[\cite{CMNR}]
If~$D(\Gamma)=(E(\Gamma);\Phi(\Gamma))$ is the delta-matroid of an embedded graph~$\Gamma$
and $E'\subset E(\Gamma)$ is a subset of its edges such that the corresponding spanning
subgraph is connected, then $D_{E'}$ is the delta-matroid of the spanning subgraph $(V(\Gamma);E')$.
Moreover, if $E'\subset E(\Gamma)$ is an arbitrary subset, and $\Gamma'_1,\dots,\Gamma'_k$
are the connected components of the corresponding spanning subgraph of~$\Gamma$, then the delta-matroid
$D(\Gamma)_{E'}$ coincides with the product of the delta-matroids $D(\Gamma'_1)\dots D(\Gamma'_k)$.
\end{proposition}

\begin{theorem}[\cite{CMNR}]
For a binary delta-matroid $D=(E;\Phi)$, its restriction $D_{E'}$ to an arbitrary
subset $E'\subset E$ is a binary delta-matroid.
\end{theorem}

The following statement shows that the coproduct defined above is compatible with the product.

\begin{proposition}
Let $D_1=(E_1;\Phi_1)$, $D_2=(E_2;\Phi_2)$ be two
delta-matroids. Then
$$
\mu(D_1D_2)=\mu(D_1)\mu(D_2).
$$
\end{proposition}

{\bf Proof.} Consider a subset $E'\subset E_1\sqcup E_2$.
Such a subset is represented as $E'=E'_1\sqcup E'_2$ with $E'_1\subset E_1$,
$E'_2\subset E_2$.
Therefore,
$$
\mu(D_1D_2)=\sum_{E'_1\subset E_1,E'_2\subset E_2}{D_1}_{E'_1}{D_2}_{E'_2}\otimes
{D_1}_{E_1\setminus E'_1}{D_2}_{E_2\setminus E'_2},
$$
since $(D_1D_2)_{E'_1\sqcup E'_2}={D_1}_{E'_1}{D_2}_{E'_2}$.
Therefore,
$$
\mu(D_1D_2)=\sum_{E'_1\subset E_1}{D_1}_{E'_1}\otimes{D_1}_{E_1\setminus E'_1}
\sum_{E'_2\subset E_2}{D_2}_{E'_2}\otimes
{D_2}_{E_2\setminus E'_2},
$$
The converse statement also is clear, which proves the Proposition.

The coproduct~$\mu$ extends by linearity to a comultiplication of the graded vector space
spanned freely by the delta-matroids. Below, we will use it only for binary delta-matroids,
and we consider the comultiplication
$$
\mu:\cB\to\cB\otimes\cB.
$$
The counit for the comultiplication is the algebra homomorphism $\cB\to\C$, which is isomorphism
when restricted to~$\cB_0$, and zero when restricted to~$\cB_i$ for $i=1,2,\dots$.

The proof of the following theorem is a routine checking of axioms, which we omit.

\begin{theorem}
The vector space~$\cB$ endowed with the comultiplication~{\rm(\ref{ecomult})}
and the multiplication~{\rm(\ref{emult})} is a graded commutative cocommutative
Hopf algebra. The subalgebra~$\cB^e\subset\cB$ spanned by even binary delta-matroids forms a Hopf
subalgebra in this Hopf algbera.
\end{theorem}

According to the Milnor--Moore theorem, each commutative cocommutative graded Hopf algebra
is nothing but the Hopf algebra of polynomials in its primitive elements. Recall that
an element~$p$ of a Hopf algebra is {\it primitive\/} if
$$
\mu(p)=1\otimes p+p\otimes 1,
$$
and that primitive elements form a vector subspace in the algebra.
The elements $s_{11},s_{12},s_{13}$ in~$\cB_1$ are primitive, and $\cB_{1}$
coincides with its primitive subspace. The elements $s_{21},s_{22},s_{23},s_{24},s_{25}$
are not primitive. Nevertheless, the dimension of the primitive subspace in~$\cB_2$ is~$5$:
any space~$\cB_n$ can be represented as the direct sum of its primitive subspace
and subspace of decomposable elements, which is spanned by polynomials in elements
of smaller degrees. In~$\cB_2$, the subspace spanned by decomposable elements is
$6$-dimensional and spanned by $s_{11}^2,s_{12}^2,s_{13}^2,s_{11}s_{12},s_{11}s_{13},s_{12}s_{13}.$

Similarly, $\cB^e_1$ coincides with its subspace of primitive elements and is $2$-dimensional,
while~$\cB^e_2$ is the direct sum of the $3$-dimensional subspace spanned by
decomposable elements and the $2$-dimensional primitive subspace.

Due to the proposition below, the Hopf algebra structure above can be restricted to binary delta-matroids such
that the empty set is feasible.

\begin{proposition}
Let $D=(E;\Phi)$ be a binary delta-matroid such that the empty set is feasible, $\emptyset\in\Phi$.
Then its restriction to any subset also is feasible.
\end{proposition}

Indeed, $D$ cannot contain coloops: otherwise $\emptyset$ would not be feasible.
And if~$e\in E$ is not a coloop, then $\emptyset$ is a feasible set for~$D\setminus e$
as well.

Therefore, both multiplication and comultiplication in~$\cB$ and~$\cB^e$
preserve the subspaces spanned by binary delta-matroids with feasible emptysets.
We denote the corresponding Hopf algebras by~$\cK=\cK_0\oplus\cK_1\oplus\cK_2\oplus\dots$ and
$\cK^e=\cK^e_0\oplus\cK^e_1\oplus\cK^e_2\oplus\dots$ (the notation reflects the fact that these
Hopf algebras are related to chord diagrams and embedded graphs with a single vertex,
that is, to knots, rather than to links). The corresponding dimensions of the spaces
of primitive elements are~$2$ for $\cK_1$, $3$ for~$\cK_2$, $1$ for~$\cK^e_1$, and $1$ for $\cK^e_2$.

\section{Four-term relations}

Vassiliev's theory of finite order knot invariants~\cite{V90} associates to a knot invariant
of order at most~$n$ a {\it weight system} of order~$n$, that is, a function on chord
diagrams (= embedded graphs with a single vertex) with~$n$ chords satisfying $4$-term relations.
This construction has a straightforward
generalization to chord diagrams of links, which are essentially embedded graphs with
the number of vertices equal to the number of connected components of the link.

The definition of the $4$-term relations requires the definition of two operations,
namely, exchanging of handle ends (the first Vassiliev move) and handle sliding (the second Vassiliev move).
The handle sliding for binary delta-matroids was defined in~\cite{MMB}. Below, we give the description
of this operation, and define the operation of exchanging handle ends.
As a result, we can introduce $4$-term relations for binary delta-matroids and the
corresponding Hopf algebra.

It was shown in~\cite{MMB} that for the delta-matroids of embedded graphs, the operation
of handle sliding, when applied to two ribbons with neighboring ends,
coincides with the handle sliding for embedded graphs. We prove a similar statement
for the operation of exchanging handle ends. Although handle sliding
and exchanging handle ends do not preserve the class of delta-matroids of embedded graphs,
they preserve a wider class of binary delta matroids. As a result, we are able to construct
a Hopf algebra of binary delta-matroids modulo $4$-term relations.

Any function
on binary delta-matroids satisfying the $4$-term relations defines a weight system,
whence a knot invariant. Therefore, studying these functions can help to construct knot invariants
and clarify their nature. We prove that any invariant of binary delta-matroids satisfying
so-called topological Tutte relations satisfies also the $4$-term relations.
This means, in particular, that the Bollobas--Riordan polynomial of delta-matroids and its relatives
produce link invariants (which was proved for a special case in~\cite{BR01}).

Note that the connected sum of chord diagrams is well defined only if
$4$-term relations are imposed. This property allows one to define the Hopf algebra
of chord diagrams modulo $4$-term relations. It was asked in~\cite{L} whether imposing
the $4$-term relations allows one to define multiplication on framed chord diagrams as well.
Recently, D.~P.~Ilyutko and V.~O.~Manturov~\cite{IM} answered this question in negative.
The results of the present section show, however, that on the level of (binary)
delta-matroids we obtain Hopf algebra structures not only for framed chord
diagrams, but for arbitrary embedded graphs as well. Multiplication in these Hopf
algebras is well defined independently of whether the $4$-term relations are imposed.

\subsection{The second Vassiliev move: handle sliding}

Let $D=(E;\Phi)$ be a set system, $a,b\in E$ be two different elements.

\begin{definition}[\cite{MMB}]
The result of {\it sliding of the element~$a$ over the element~$b$\/}
is the set system $\widetilde D_{ab}=(E;\widetilde\Phi_{ab})$, where
$\widetilde\Phi_{ab}=\Phi\Delta\{\phi\sqcup\{a\}|\phi\sqcup\{b\}\in\Phi
\text{ and } \phi\subset E\setminus\{a,b\}\}$.
\end{definition}

It is proved in~\cite{MMB} that if $D=(E(\Gamma);\Phi(\Gamma))$ is the delta-matroid
of an embedded graph~$\Gamma$ and $a,b$ are two ribbons in~$\Gamma$ with neighboring ends,
then the delta-matroid of the ribbon graph $\widetilde\Gamma_{ab}$
obtained from~$\Gamma$ by sliding the handle~$a$ over the handle~$b$
coincides with the delta-matroid $\widetilde D_{ab}$. However, if
the ends of the ribbons $a,b$ in~$\Gamma$ are not neighboring, then the handle
sliding of the above definition can lead to a set system that is not
isomorphic to the delta-matroid of any embedded graph.
Moreover, the following example from~\cite{MMB} shows that a handle sliding
applied to a delta-matroid can produce a set system that is not a delta-matroid.

\begin{example}
For the delta-matroid $$D=(\{1,2,3\};\{\emptyset,\{1,2\},\{1,3\},\{2,3\},\{1,2,3\}\})$$
the set system $\widetilde D_{12}=(\{1,2,3\};\{\emptyset,\{1,2\},\{2,3\},\{1,2,3\}\})$
is a delta-matroid no longer.
\end{example}

Nevertheless, the following
theorem is valid.

\begin{theorem}[\cite{MMB}]\label{tV2}
If $D=(E;\Phi)$ is a binary delta-matroid and $a,b$ are two distinct elements in~$E$,
then $\widetilde D_{ab}$ is a binary delta-matroid.
\end{theorem}

In the next section we prove a similar theorem for the other Vassiliev move, the
first one.

In~\cite{L00}, the second Vassiliev move was defined for abstract graphs.
We are going to show that this definition is, in fact, consistent with the definition above.
Let us recall the definition from~\cite{L00} (together with its extension to framed graphs in~\cite{L}).
For a framed abstract graph~$G$ and a pair of vertices~$a,b\in V(G)$ in it, the
graph~$\widetilde G_{ab}$ is defined as a graph on the same set $V(G)$ of vertices
such that the adjacency of any vertex~$c$ to~$a$, $c\ne a,b,$ toggles iff
$c$ is adjacent to~$b$ in~$G$.  In addition, the adjacency of~$a$ and~$b$ toggles
if the framing of~$b$ is~$1$.

If~$G$ is the intersection graph of a chord diagram, and $a,b$ are two chords with neighboring
ends in the diagram, then this move indeed corresponds to sliding of the handle~$a$
along the handle~$b$~\cite{L}.

\begin{theorem}
For an abstract framed graph~$G$, we have
$$
D(\widetilde G_{ab})=\widetilde {D(G)}_{ab}.
$$
\end{theorem}

{\bf Proof.} Indeed, the adjacency matrix $A(G)$ of an abstract framed graph~$G$
can be considered as the matrix of a symmetric binary form over the field of
two elements~$\F_2$ in the vector space $\F_2^{V(G)}$ spanned by the vertices of the graph.
The second Vassiliev move $G\mapsto\widetilde G_{ab}$ does not modify the form,
but changes the basis:
$$
(a,b,c,\dots)\mapsto (a+b,b,c,\dots).
$$
(Note that this property justifies the name of the move: on the homology level the
second Kirby move in topology of $3$-manifolds does exactly the same thing, but over~$\Z$
rather than over~$\F_2$).

Of course, this change of basis does not affect the (non)degeneracy
property of any subset of vertices in~$G$ not containing~$a$
or containing both $\{a,b\}$. Now, if a subset $U\sqcup\{b\}\subset V(G)$
does not contain~$a$ and is nondegenerate, then the nondegenercy of $U\sqcup\{a\}$
toggles between $G$ and $\widetilde G_{ab}$.

\subsection{The first Vassiliev move: exchanging handle ends}

For an embedded graph~$\Gamma$ and two distinct ribbons $a,b\in E(\Gamma)$
such that one of the ends of~$a$ is a neighbor of one of the ends of~$b$ along some
vertex, the first Vassiliev move consists in exchanging these neighboring ends.
The following definition mimics what happens with the underlying delta-matroids
under this operation.

Let $D=(E;\Phi)$ be a set system, $a,b\in E$ be two different elements.

\begin{definition}
The result of {\it exchanging of the ends of the ribbon~$a$ and the ribbon~$b$\/}
is the set system $D'_{ab}=(E;\Phi'_{ab})$, where
$\Phi'_{ab}=\widetilde{(\Phi*b)}_{ab}*b$.
\end{definition}

Note that, in contrast to the second Vassiliev move, the first Vassiliev move is
symmetric with respect to the ribbons~$a,b$ whose neighboring ends we exchange,
$D'_{ab}=D'_{ba}$.

Since the operation~$*$ preserves the class of binary delta-matroids, Theorem~\ref{tV2}
immediately implies

\begin{proposition}
If $D=(E;\Phi)$ is a binary delta-matroid and $a,b$ are two distinct elements in~$E$,
then $D'_{ab}$ is a binary delta-matroid.
\end{proposition}

\begin{theorem}\label{tfVmeg}
If $D=(E(\Gamma);\Phi(\Gamma))$ is the delta-matroid
of an embedded graph~$\Gamma$ and $a,b$ are two ribbons in~$\Gamma$ with neighboring ends,
then the delta-matroid of the ribbon graph $\Gamma'_{ab}$
obtained from~$\Gamma$ by exchanging the ends of the handles~$a$ and~$b$
coincides with the delta-matroid $D'_{a,b}$.
\end{theorem}

{\bf Proof.} The set system $D*b$ is the delta-matroid of the partial dual embedded
graph $\Gamma^b$, see~\cite{C09} or~\cite{CMNR}. After taking the partial dual along~$b$,
sliding the neighboring end of the handle~$a$ along the new~$b$ and returning~$b$
to its original place, we obtain exactly the neighboring ends exchange move.

Vassiliev moves for binary delta-matroids possess properties similar to those
for embedded graphs:

\begin{proposition}
The following statements about the Vassiliev moves are valid:
\begin{itemize}
\item the first Vassiliev move is an involution, $(D'_{ab})'_{ab}=D$;
\item the second Vassiliev move is an involution, $\widetilde{(\widetilde D_{ab})}_{ab}=D$;
\item the first and the second Vassiliev moves commute, $(\widetilde D_{ab})'_{ab}=\widetilde{(D'_{ab})}_{ab}$.
\end{itemize}
\end{proposition}

\subsection{The four-term relation for binary delta-matroids}

As usual, we say that an invariant~$f$ of embedded graphs {\it satisfies the four-term relation\/}
if for any embedded graph~$\Gamma$ and any pair $a,b$ of its distinct edges having neighboring
ends we have
\begin{equation}\label{e4term}
f(\Gamma)-f(\Gamma'_{ab})=f(\widetilde\Gamma_{ab})-f(\widetilde\Gamma'_{ab}).
\end{equation}

Similarly, we say that an invariant~$f$ of binary delta-matroids {\it satisfies the four-term relation\/}
if for any binary delta-matroid~$D$ and a pair of distinct elements~$a,b$ in its ground set we have
\begin{equation}\label{e4term}
f(D)-f(D'_{ab})=f(\widetilde D_{ab})-f(\widetilde D'_{ab}).
\end{equation}

Theorem in~\cite{MMB} and Theorem~\ref{tfVmeg} above mean that

\begin{theorem}
Any invariant of binary delta-matroids satisfying the $4$-term relation~{\rm(\ref{e4term})}
defines a weight system, whence a link invariant.
\end{theorem}

In the next section we show that, particularly, each Tutte invariant of delta-matroids
defines a link invariant.

\subsection{Tutte relations for invariants of binary delta-matroids and weight systems}

We say that an invariant~$f$ of delta-matroids {\it satisfies the Tutte relations\/}
if, for any delta-matroid $D=(E;\Phi),$  we have
\begin{eqnarray}\label{eTutte}
f(D)&=&x f(D\setminus e)+ y f(D/e) \text{ for any }e\in E\nonumber\\
&&\text{ that is neither loop nor coloop;}\nonumber\\
f(D)&=&z f(D\setminus e)\text{ for any loop }e\in E;\\
f(D)&=&w f(D/e) \text{ for any coloop }e\in E.\nonumber
\end{eqnarray}
Here $x,y,z,w$ are some indeterminates. Note that, since a minor of a binary delta-matroid
also is a binary delta-matroid, we may as well consider invariants of binary delta-matroids
satisfying the Tutte relations.

A generic example of an invariant satisfying the Tutte relations is the Tutte polynomial
of delta-matroids.
The Tutte relation for graphs on surfaces appeared first in paper~\cite{BR01}, where
the Tutte polynomial for embedded graphs on orientable surfaces has been introduced.
It was proved in~\cite{BR01} that the Tutte polynomial satisfies the $4$-term relation
for orientable embedded graphs and generates therefore a knot invariant.
In this section we prove a generalization of this statement showing that any invariant
of binary delta-matroids satisfying the Tutte relation satisfies also the $4$-term relation.
This means, in particular, that the Tutte polynomial for delta-matroids satisfies
the $4$-term relation and defines thus a link invariant. Since the precise definition
of the Tutte polynomial requires some additional notions not considered in the present
text, we refer the reader to~\cite{CMNR} for it.

\begin{proposition}
If~$f$ is an invariant of binary delta-matroids satisfying the Tutte relation~{\rm(\ref{eTutte})},
then~$f$ satisfies the $4$-term relation~{\rm(\ref{e4term})}.
\end{proposition}

{\bf Proof.} Let $D=(E;\Phi)$ be a binary delta-matroid.
The proof of the proposition requires the following two relations.
If $a,b\in E$ are arbitrary distinct elements, then
\begin{itemize}
\item $D'_{ab}\setminus b=D\setminus b$;
\item $\widetilde D_{ab}/b=D/b$.
\end{itemize}
Both are obvious.

Now let $f$ be an invariant of binary delta-matroids satisfying the Tutte relations.
We are going to prove that
$$
f(D)-f(D'_{ab})=f(\widetilde D_{ab})-f(\widetilde D'_{ab})
$$
for arbitrary pair of distinct elements $a,b\in E$.

If~$b$ is neither a loop, nor a coloop in~$D$, then
it is neither a loop, nor a coloop in~$D'_{ab}, \widetilde D_{ab}, \widetilde D'_{ab}$,
and we have
$$
f(D)=x f(D\setminus b)+y f(D/b),
$$
and similarly for the other three terms of the $4$-term relation.
Make this substitution and take first the coefficients of~$x$.
Due to the equation $D'_{ab}\setminus b=D\setminus b$, both the left- and
the right-hand side of the equation are~$0$.
In its turn, the coefficient of~$y$ on the left-hand side has the form
$f(D/b)-f(D'_{ab}/b)$ and it coincides with the coefficient of~$y$
on the right-hand side due to the equation $\widetilde D_{ab}/b=D/b$.

The other cases are considered in a similar way.

Thus, we can prove the proposition by induction on the number of elements in the
ground set of the binary delta-matroid.

\subsection{Hopf algebras of binary delta-matroids modulo $4$-term relations}

The Hopf algebra~$\cB$ of binary delta-matroids, as well as its Hopf subalgebra $\cB^e$
of even binary delta-matroids can be factorized modulo the $4$-term relations.
Denote by $\cF\cB$ (respectively, $\cF\cB^e$) the graded quotient space
of the space of binary matroids (respectively, even binary matroids) modulo the
$4$-term relations:
\begin{eqnarray*}
\cF\cB_i&=&\cB_i/\langle  D-D'_{ab}-\widetilde D_{ab}+\widetilde D'_{ab}\rangle, \quad i=0,1,2,\dots\\
\cF\cB^e_i&=&\cB^e_i/\langle  D-D'_{ab}-\widetilde D_{ab}+\widetilde D'_{ab}\rangle,
\quad i=0,1,2,\dots.
\end{eqnarray*}

\begin{theorem}
The multiplication~$m$ and the comultiplication~$\mu$ induce on the spaces $\cF\cB$
and $\cF\cB^e$ the structure of graded commutative cocommutative Hopf algebras.
\end{theorem}

\begin{example}
The vector spaces $\cF\cB^e_i$ for $i=0,1,$ and~$2$ coincide with the vector space~$\cB^e_i$,
since the even $4$-term relations are trivial for these values of~$i$.
In contrast, there is a (single) nontrivial $4$-term relation for $i=2$ in the noneven case:
$$
s_{11}s_{12}-s_{22}=s_{23}-s_{12}^2.
$$
Therefore, $\cF\cB_2=\cB_2/\langle s_{11}s_{12}-s_{22}-s_{23}+s_{12}^2\rangle$, $\dim \cF\cB_2= 10$, and the primitive subspace in it is $4$-dimensional.
Indeed, none of the elements $s_{22},s_{23}$ is decomposable, but their sum is.
\end{example}

Since both the first and the second Vassiliev move preserve the class of binary delta-matroids with feasible empty set,
the quotients $\cF\cK$ and $\cF\cK^e$ of the Hopf algebras $\cK$ and $\cK^e$, respectively, modulo
the $4$-term relations also are Hopf algebras. For $n=1,2$ the corresponding $4$-term relations are trivial.

Let us collect the computed dimensions of the spaces of primitive elements into a table.
\begin{table}[h]
\begin{center}
\begin{tabular}{c|c|c}
n&1&2\\
\hline $\cB_n$&3&5\\
\hline $\cB^e_n$&2&2\\
\hline $\cF\cB_n$&3&4\\
\hline $\cF\cB^e_n$&2&2\\
\hline $\cK_n$&2&3\\
\hline $\cK^e_n$&1&1\\
\hline $\cF\cK_n$&2&3\\
\hline $\cF\cK^e_n$&1&1\\
\end{tabular}
\end{center}
\caption{Dimensions of the primitive subspaces}
\end{table}

\begin{example}
The weight system~$w_C$ on framed chord diagrams corresponding to the Conway invariant of knots
can be defined as the function taking on a chord diagram value~$1$ if the corresponding
one-vertex ribbon graph has a connected boundary and~$0$ otherwise. This weight system
admits a natural extension to binary delta-matroids: for a binary delta-matroid
$D=(E,\Phi)$, define $w_C(D)=1$ if $E\in\Phi$ and~$0$ otherwise. This function
satisfies the $2$-term relation: $w_C(D)=w_C(\widetilde D_{ab})$ for any pair of distinct
elements $a,b\in E$, whence the $4$-term relation. We extend it to~$\cF\cB$
by linearity.

The function~$w_C$ obviously is multiplicative, $w_C(D_1,D_2)=w_C(D_1)w_C(D_2)$
for any pair of binary delta-matroids $D_1,D_2$. Therefore, its logarithm 
is well defined. The value of this logarithm on chord diagrams is known to 
be related to the weight system $\sl_2$, see details in~\cite{BNV,KLMR}.
Hence, the value of $\log w_C$ on binary delta-matroids can be considered
as a manifestation of the existence of a yet unknown construction
of an $\sl_2$-weight system on binary delta-matroids extending that for chord diagrams. 
This construction is unknown yet even for (framed) graphs, see~\cite{KLMR}.
\end{example}

\subsection{Vassiliev moves and Lagrangian subspaces}

In~\cite{KS} it is shown that the first and the second Vassiliev moves
for embedded graphs can be naturally expressed as
base changes in the $2|E|$-dimensional symplectic space over~$\F_2$ spanned by the edges of the
graph and their duals. We reproduce the definition of these base changes below and show that
it is compatible with the above definition of the Vassiliev moves for
binary delta-matroids.

Let $D=(E;\Phi)$ be a delta-matroid. Denote by~$E^*$ a copy of~$E$;
the element of~$E^*$ that is a copy of an element~$e\in E$ is denoted
by~$e^*$. Denote by~$V_E$ the $2|E|$-dimensional vector space over~$\F_2$
spanned by~$E\sqcup E^*$. Introduce a bilinear form $(\cdot,\cdot)$ on~$V_E$ by the rule
$(e,e^*)=(e^*,e)=1$ for any $e\in E$, all the other pairings between basic vectors being zero.
This form is nondegenerate and skew-symmetric, thus making~$V_E$ into
a symplectic space.

\begin{definition}[\cite{KS}]
For $a,b\in E$, $a\ne b$, the {\it first Vassiliev move\/} $T_1^{ab}:V_E\to V_E$ of the
space~$V_E$ is defined by
$$
a\mapsto a,\quad b\mapsto b,\quad
a^*\mapsto a^*+b,\quad
b^*\mapsto b^*+a,
$$
being identical on the other basic vectors.

For $a,b\in E$, $a\ne b$, the {\it second Vassiliev move\/} $T_2^{ab}:V_E\to V_E$ of the
space~$V_E$ is defined by
$$
a\mapsto a+b,\quad b\mapsto b,\quad
a^*\mapsto a^*,\quad
b^*\mapsto a^*+b^*,
$$
being identical on the other basic vectors.
\end{definition}

Now, to each subset $E'\subset E$ the coordinate subspace $L_{E'}\subset V_E$
can be associated; this subspace is spanned by the basic vectors
corresponding to the elements in $E'$ as well as the elements in $E^*\setminus (E')^*$.
Each such subspace is $|E|$-dimensional and isotropic, meaning that $(v_1,v_2)=0$
for any pair of vectors $v_1,v_2\in L_{E'}$, whence a Lagrangian subspace in~$V_E$.

\begin{proposition}
Let $D=(E,\Phi)$ be a binary delta-matroid, $a,b\in E$, $a\ne b$.
Then the actions of the operations $T_1^{ab},T_2^{ab}$ on the symplectic
vector space $V_E$ induce on~$D$ the first and the second Vassiliev move,
respectively:
$$
T_1^{ab}:D \mapsto D'_{ab},\qquad T_2^{ab}:D \mapsto \widetilde D_{ab}.
$$
\end{proposition}

The proof is straightforward.

In~\cite{KS} the action of Vassiliev moves on Lagrangian subspaces in~$V_E$
was used to introduce the $4$-term relations and the Hopf algebra of Lagrangian
subspaces modulo the $4$-term relations. There is a mapping taking an embedded graph to a Lagrangian
subspace, and linear functionals on this Hopf algebra
determine weight systems.

In our construction, the same Vassiliev moves act on a tuple of Lagrangian subspaces
corresponding to the feasible sets of a binary delta-matroid rather than on a single
Lagrangian subspace.


\begin{thebibliography}{99}

\bibitem{BNV} D.~Bar-Natan, H.~Vo, {\it Proof of a conjecture of Kulakova et al. 
related to the $\sl_22$ weight system}, European Journal of Combinatorics,
Volume 45, April 2015, Pages 65--70, arXiv:1401.0754

\bibitem{BR01} B.~Bollob\'as, O.~Riordan,
{\it A polynomial invariant of graphs on orientable surfaces},
Proc. London Math. Soc. (3) 83 (2001), no. 3, 513–531.

\bibitem{BR02} B.~Bollob\'as, O.~Riordan,
A polynomial of graphs on surfaces.
Math. Ann. 323 (2002), no. 1, 81–96.

\bibitem{Bo871} A.Bouchet,
{\it Greedy algorithms and symmetric matroids},
Math. Programm. 38 (1987), 147--159

\bibitem{Bo872} A.Bouchet,
{\it Representability of $\Delta$-matroids},
in: Proceedings of the 6th Hungarian Colloquium of Combinatorics,
Colloq. Math. Soc. J\'anos Bolyai 38 (1987), 167--182

\bibitem{Bo89} A.Bouchet,
{\it Maps and delta-matroids},
Discrete Math. 78 (1989), 59--71


\bibitem{C09} S.~Chmutov,
{\it Generalized duality for graphs on surfaces and the signed Bollob\'as--Riordan
polynomial},
J. of Combin. Theory Ser. B 99 (2009) 617--638


\bibitem{CL07} S.~Chmutov, S.~Lando,
{\it Mutant knots and intersection graphs},
Algebraic \& Geometric Topology 7 (2007) 1579–-1598

\bibitem{CMNR} C.~Chun, I.~Moffatt, S.~D.~Noble, R.~Rueckriemen,
{\it Matroids, delta-matroids and embedded graphs},
arXiv: 1403.0920v1, 45 pp.

\bibitem{IM} D.~P.~Ilyutko, V.~O.~Manturov,
{\it A parity map of framed chord diagrams},
arXiv: 1506.0918

\bibitem{KS} V.~Kleptsyn, E.~Smirnov
{\it Ribbon graphs and bialgebra of Lagrangian subspaces}, arxiv:1401.6160

\bibitem{KLMR}     E.~Kulakova, S.~Lando, T.~Mukhutdinova, G.~Rybnikov, 
{\it On a weight system conjecturally related to $\sl_2$},
European Journal of Combinatorics, Volume 41, October 2014, Pages 266--277 


\bibitem{L00}      {S.~K.~Lando,}
    	{\it On a Hopf algebra in graph theory},
    	J. Comb. Theory, Ser. B, vol.~{\bf 80} (2000), 104--121.

\bibitem{L}      {S.~K.~Lando,}
    	{\it $J$-invariants of ornaments and framed chord diagrams},
    	Funct. Anal. Appl.,~{\bf 40}(1) (2006), 1--13.


\bibitem{LZ04} S.~Lando, A.~Zvonkin, {\it Graphs on surfaces and their applications}, Springer, 2004.

\bibitem{MMB} Iain Moffatt, Eunice Mphako-Banda,
{\it Handle slides for delta-matroids}, arXiv:1510.07224, 12 pp.

\bibitem{S01} E.~Soboleva
{\it Vassiliev knot invariants coming from Lie algebras and 4-invariants},
Journal of Knot Theory and Its Ramifications, Volume 10, Issue 01, February 2001,
161--169


\bibitem{V90}
V.~A.~Vassiliev,
{\it Cohomology of knot spaces}, in: Theory of singularities and its applications, 23–-69,
Adv. Soviet Math., 1, Amer. Math. Soc., Providence, RI, 1990.

\end{thebibliography}
\end{document}